\documentclass{amsart}
\usepackage{amssymb,latexsym}

\newtheorem{theorem}{Theorem}[section]
\newtheorem{lemma}[theorem]{Lemma}
\newtheorem{proposition}[theorem]{Proposition}
\newtheorem{cor}[theorem]{Corollary}

\theoremstyle{definition}
\newtheorem{definition}[theorem]{Definition}

\theoremstyle{remark}
\newtheorem{remark}[theorem]{Remark}

\numberwithin{equation}{section}

\newcommand{\soit}{\stackrel{def}{=}}
\newcommand{\Hom}{\mbox{Hom}\,}
\newcommand{\Ext}{\mbox{Ext}\,}

\newcommand{\Tr}{\textup{Tr}}
\newcommand{\wt}{\textup{wt}}
\newcommand{\cS}{\mathcal {S}}
\newcommand{\sop}{\mathit{Op}}
\newcommand{\sopr}{\mathit{Opr}}
\newcommand{\sopc}[1]{\sop(#1)}

\newcommand{\za}{\alpha}

\newcommand{\zO}{\Omega}

\newcommand{\zS}{\Upsilon}

\DeclareSymbolFont{lasy}{U}{lasy}{m}{n}
\SetSymbolFont{lasy}{bold}{U}{lasy}{b}{n}
\let\Box\undefined 
\DeclareMathSymbol\Box {\mathord}{lasy}{"32}

\newcommand{\suitr}[3]{0\to {#1}\to {#2}\to{#3}\to 0}

\newcommand{\orb}[1]{{\mathcal O}_\mathbf{ #1}}

\begin{document}

\title[Rational smoothness]{Rational smoothness of
varieties of representations for quivers of Dynkin type}

\address{}
\curraddr{}
\email{}
\thanks{}


\author{Philippe Caldero}
\address{D\'epartement de Math\'ematiques, Universit\'e Claude Bernard Lyon
I,
69622 Villeurbanne Cedex, France}
\email{caldero@igd.univ-lyon1.fr}


\author{Ralf Schiffler}
\address{School of mathematics and statistics \\ Carleton University\\
1125 Colonel By Drive\\ Room 4302 Herzberg Building\\
Ottawa, Ontario\\ Canada K1S 5B6}
\email{ralf@math.carleton.ca}
\thanks{The second author was supported in part by FCAR Grant.}

\begin{abstract}
Let $\mathbf {U}^+$ be the positive part of the quantized
enveloping algebra $\mathbf{U}$
of type $A,D$ or $E$.
The change of basis between canonical, and PBW-basis of
$\mathbf{ U}^+$ has a geometric interpretation in terms of local intersection
co\-ho\-mo\-logy of some affine algebraic varieties, namely
the Zariski closures of orbits
of representations of a quiver of type $A,D$ or $E$.
In this paper we characterize
the rationally smooth orbit closures and prove in particular that
orbit closures are smooth if and only if they are
rationally smooth. This provides an analogue of theorems of V. Deodhar, and J. Carrell and D. Peterson on Schubert varieties.
\end{abstract}

\maketitle

\begin{section} {Introduction}
Let $\mathbb{F}$ be an algebraically closed field,
 $\mathbf{ d} = (d_1, d_2, \dots, d_n) \in \mathbb{N}^n$
and $G_\mathbf{ d} = \prod_{i = 1}^n GL_{d_i}(\mathbb{F})$.
Let
$\mathcal {Q}$ be a fixed quiver whose underlying graph $\Delta$ is the
Dynkin graph of type $A_n,D_n$ or $E_n$.
$G_\mathbf{ d}$ acts on $E_\mathbf{ d} = \oplus_{i\rightarrow j \in
\mathcal{ Q}} \Hom_\mathbb{F}(\mathbb{F}^{d_i}, \mathbb{F}^{d_j})$, by conjugation.
Let $\mathcal{O}$ be a $G_\mathbf{ d}$-orbit and $\overline{\mathcal{O} }$ its
Zariski closure.
In \cite{bedschi1} the complete list of rationally smooth orbit
closures of type $A_n$ was obtained. As a consequence it was shown that, in
type $A_n$,
$\overline{\mathcal{O} }$ is
rationally smooth if and only if $\overline{\mathcal{O} }$ is smooth.
In this paper we will generalize these results to the types $D_n$ and
$E_n$, see theorem \ref{smoothness} and corollary \ref{cor}.

Rational smoothness is a
topological property, which is defined using local
intersection cohomology, and has been extensively studied for
Schubert varieties, \cite{deodhar1}, \cite{carrell1}. For a survey of some of these results,
see \cite{BL} and \cite{brion1}.

Let $\mathbf{ U}^+$ be the positive part of the quantized
enveloping algebra $\mathbf{ U}$ over $\mathbb{ Q}(v)$ associated
by Drinfeld and Jimbo to the root system of type $\Delta$.
Kashiwara and Lusztig have constructed independently of each
other a unique canonical basis $\mathbf{ B}$ of $\mathbf{ U}^+$
in \cite{kashiwara1} and \cite{lusztig1}.
For each reduced expression $\mathbf{ i}$ of the
longest element
$w_0$ of the Weyl group $W$ of type
$\Delta$, there is also a PBW-basis $B_\mathbf{ i}$. Some of the
reduced expressions are adapted to the quiver $\mathcal {Q}$. In
this case, Lusztig has shown in \cite{lusztig1} that the entries of
the transition matrix between the bases $\mathbf{ B}$ and
$B_\mathbf{ i}$ have a description in terms of local intersection
cohomology of orbit closures.
We use this approach to study rational smoothness of orbit closures.
One important ingredient therefore is the action of the bar involution of
$\mathbf{U}$ on the elements of the PBW-basis. 

This paper is organized as follows.
In section 2 we fix notations and recall some results that we will need at
a later stage. In particular, we recall the Hall algebra realization of $\mathbf{ U}^+$.
We then present two different approaches to study rational smoothness of
orbit closures: one algebraic and the other geometric.\par\noindent
Section 3 contains the algebraic approach. Here we consider
the bar of a PBW-basis element as a linear combination in the
PBW-basis, and study the coefficients of this expansion; to be more
precise, we
calculate the derivative at $v=1$ of these coefficients.
The method used to calculate the coefficients differs from \cite{bedschi1}.
In fact, we use the dual approach, with the help of the canonical form on
$\mathbf{ U}^+$. Indeed, the dual PBW-basis is known in the Hall algebra
context by \cite{green1}. Moreover, the adjoint $\sigma$ of the bar
involution is a $\mathbb{Q}$-antiautomorphism (up to a power of $q$) given
by Lusztig, \cite{lusztig3}.
The reason why things are more convenient when we dualize is that the dual
of an irreducible element of the PBW-basis is an element of the dual
canonical basis, \cite{caldero1}, and so is stable by $\sigma$. Then, the
image by $\sigma$ of a general dual PBW-basis element can be easily
calculated by the antiautomorphism property. This enables us to realize the
desired coefficients in terms of generalized Hall polynomials, see
Proposition \ref{formula}. Using the same methods as in \cite{bedschi1},
we then obtain the complete list of rationally smooth orbit closures. It is easy
to see that each orbit closure in this list is smooth.

\par\noindent
In section 4 we present another proof of this characterization using a
geometric approach. According to an idea of Michel Brion, we calculate the
Euler-Poincar\'e characteristic of the projectivization of the orbit
closures $\overline{\mathcal{O} }$. Indeed, it is known that, for
rationally smooth cones, the Euler-Poincar\'e characteristic of the
projectivization of the cone equals the dimension of the cone. Thanks to a
theorem of Deligne on the Weil conjecture, the characteristic can be
calculated by counting the number of $\mathbb{F}_{q}$-rational points of an
orbit closure viewed as a variety on
an algebraic closure $\overline{\mathbb{F}_q}$ of $\mathbb{F}_q$,
 and specializing at
$q=1$. This method provides a geometric interpretation of the proof given
in \cite{bedschi1}, and followed in 3. It also provides an interesting
homological realization of algebraic elements, namely the derivative at
$q=1$ of the coefficients of the bar involution in the
PBW-basis.\smallskip\noindent

Our main result, Theorem \ref{main}, generalizes \cite[Theorem
5.4]{bedschi1}. Nevertheless, the E-case slightly differs from the A-D
case. Let's explain how. First of all, by Ringel's Hall algebra approach of
quantum groups, the coefficients of the bar automorphism discussed above
are parameterized by couples $(\mathcal{O},\mathcal{O'})$ of
$G_{\mathbf{d}}$ orbits, and are non zero if and only if
$\mathcal{O'}\preceq\mathcal{O}$ for the so-called degeneration ordering.
Now, Klaus Bongartz gave, \cite{bongartz1}, a description of this ordering
in the representation theory of the quiver $\mathcal{Q}$. In particular, we
can define ``elementary degenerations'' corresponding to certain non-split
exact sequences, called elementary operations, in the Auslander-Reiten
quiver of $\mathcal{Q}$, and
$\mathcal{O'}\preceq\mathcal{O}$ if and only if there exists a chain of
elementary degenerations relying $\mathcal{O'}$ and $\mathcal{O}$. The
results in \cite{bedschi1} proves that the coefficient parameterized by
$(\mathcal{O},\mathcal{O'})$ "sees" the minimal length of such a chain, and
that the derivative at $q=1$ of the coefficient is non zero if and only if
the degeneration corresponding to $(\mathcal{O},\mathcal{O'})$ is
elementary. In order to generalize this result, we have to
replace elementary degenerations by a finer notion, see 3.2.
\end{section}

\begin{section}{Notations and recollections}
\begin{subsection}{The quantized enveloping algebra $\mathbf{ U}$}
Let $\Delta$ be a Dynkin diagram of type $A,D$ or $E$, let $n$ be the number
of
vertices of $\Delta$ and let $(a_{ij})$ be the corresponding Cartan matrix.
Thus
$$a_{ij}=\left\{\begin{array}{ll}
2 & \textup{if } i=j\\
-1 & \textup{if } i{\quad \over \quad}j \textup{ is an edge in } \Delta \\
0 & \textup{otherwise.}
\end{array}\right.
$$
Let $v$ be an indeterminate and $\mathbf{ U}$ the quantized enveloping algebra
of Drinfeld-Jimbo
of type $\Delta$ over the field $\mathbb{ Q}(v)$ of rational functions.
$\mathbf{ U}$ is a $\mathbb{ Q}(v)$-algebra with generators:
$E_i,\,F_i,\,K_i,\,K_i^{-1} (1\le i\le n)$ and relations:
$$ K_i K_i^{-1} = K_i^{-1} K_i = 1,\quad K_i K_j = K_j
K_i;
\leqno\hbox{(r.1)}$$
$$K_i E_j =v^{a_{ij}} E_j K_i
\leqno\hbox{(r.2)}$$
$$K_i F_j =v^{-a_{ij}}
\leqno\hbox{(r.3)}$$
$$E_i F_j - F_j E_i = \delta_{ij} {(K_i - K_i^{-1})\over(v -v^{-1})}
\quad \textup{ where }
\delta_{ij} = \left\{\begin{array}{ll}
1 & \textup{if }i = j\\
0& \textup{if } i \ne j.
\end{array} \right.
\leqno\hbox{(r.4)}$$
$$\left\{\begin{array}{ll}
E_i^2 E_j - (v + v^{-1}) E_i E_j E_i + E_j E_i^2= 0
& \textup{if } a_{ij}=-1\\
E_i E_j - E_j E_i = 0 & \textup{if } a_{ij}\ne -1.
\end{array}\right.
\leqno\hbox{(r.5)}$$
$$\left\{\begin{array}{ll}
F_i^2 F_j - (v + v^{-1}) F_i F_j F_i + F_j F_i^2= 0
& \textup{if } a_{ij}=-1\\
F_i F_j - F_j F_i = 0,& \textup{if } a_{ij}\ne-1
\end{array}\right.
\leqno\hbox{(r.6)}$$

Let $\mathbf{ U^+}$ be the $\mathbb{ Q}(v)$-subalgebra generated by the
$E_i\quad(1\le i\le n)$.
Let $\overline{(\ )}:\mathbf{ U} \rightarrow \mathbf{ U}$ be the
involution of $\mathbb{ Q}$-algebras defined by
$$E_i \mapsto E_i, \quad F_i \mapsto F_i, \quad
K_i \mapsto K_i^{-1} \quad \hbox{ for all } 1 \leq i \leq n
\quad \hbox{ and } \quad v \mapsto v^{-1}.$$
Note that $ \overline{\mathbf{ U^+ }} =\mathbf{ U}^+$.

Let $Q$, resp. $Q^+$, be the free abelian group, resp. semigroup, with
basis
$\{\alpha_1, \alpha_2, \ldots, \alpha_n\}$.
Define an inner product
$(\ ,\ )_Q$ on $Q$ by $(\alpha_i, \alpha_j)_Q = a_{ij}$.
Let $R= \{\alpha \in Q \mid (\alpha, \alpha)_Q = 2\}$. $R$ is
a root system of type $\Delta$ whose set of simple roots
is
$\{\alpha_1, \alpha_2, \ldots, \alpha_n\}$.
Let $R^+ = \{\alpha \in R \mid \alpha = \sum_j c_j \alpha_j \hbox{ with }
c_j \in \mathbb{ N}\}$ be the subset of positive roots.

Each $\alpha \in R$ defines a reflection $s_{\alpha}:Q
\rightarrow Q$, $\ z \mapsto z - (z, \alpha)_Q\ \alpha$. We will write
$s_i$ instead of $s_{\alpha_i}$. Let $W$ be the Weyl group of $R$.
This is the subgroup of ${\rm Aut}(Q)$
generated by the reflections $s_i$, ($1 \leq i \leq n$).
Let $\ell(w)$ be the length of $w$ with respect to the
generators
$\{s_1, s_2, \ldots, s_n\}$ and denote by $w_0$ the unique element of $W$ of
maximal length. It is known that
$\ell(w_0) = \nu =\#(R^+)$.
Let $\Tr$ be the linear form on $Q^+$ such that
$\Tr(\alpha_i)=1$, $1\leq i\leq n$.

We shall use the $Q^+$-grading $\wt$ of $\mathbf{U}^+$ defined by
$\wt(E_i)=\alpha_i$.

Lusztig has defined an action of the braid group
on $\mathbf{ U}$ \cite{lusztig1} and used it to define bases of
PBW type of $\mathbf{ U}^+$. We now recall these
definitions.

For $i \in \{1,\ldots, n\}$, let $\tilde T_i:\mathbf{ U} \rightarrow \mathbf{
U}$ be the automorphism of $\mathbb{ Q}(v)$-algebras defined by

\begin{tabular}{llll}
$E_i \mapsto - K_i^{-1} F_i$ , &
$F_i \mapsto - E_i K_i$ , & $K_i \mapsto K_i^{-1}$ ,\\
$E_j \mapsto E_j$, & $F_j \mapsto F_j$, & $K_j \mapsto
K_j$, & if $a_{ij}=0$,\\
$E_j \mapsto (E_j E_i - v^{-1}E_i E_j)$,&$F_j
\mapsto (F_i F_j - v F_j F_i)$,& $K_j \mapsto K_i K_j$ , &if
$a_{ij} = -1$.\\
\\
\end{tabular}
We have $\tilde T_i \tilde T_j \tilde T_i = \tilde T_j
\tilde T_i \tilde T_j$\quad if $a_{ij}=-1$ and
$\tilde T_i \tilde T_j = \tilde T_j \tilde T_i$\quad if
$a_{ij}\ne -1$.
This gives us a braid group action.
Moreover
$\tilde T_i(E_j)
= \tilde T_j^{-1}(E_i)$ if $a_{ij}=-1$.

Given integers $M, N \geq 0$, we define
$$
[N]! = \displaystyle{\prod_{h = 1}^N {(v^h - v^{-h})\over (v - v^{-1})}}
\ \in \ \mathbb{ Z}[v, v^{-1}],
\
\left[\displaystyle{{M + N\atop N}}\right] =
\displaystyle{{[M + N]!\over [M]![N]!}}\ \in \ \mathbb{ Z}[v, v^{-1}]$$
$$ \hbox{ and }\qquad E_i^{(N)} = \displaystyle{{E_i^N\over [N]!}}
\hbox{ for } 1
\leq i \leq n.
$$

Let ${\mathcal I}$ be the set of sequences
$\mathbf{ i} = (i_1, \ldots, i_{\nu})$ of elements in
$ \{1,\ldots, n\}$ such that
$s_{i_1} \ldots s_{i_{\nu}}$ is a reduced expression of $w_0$.
Each $\mathbf{ i} \in {\mathcal I}$ gives rise to a total order on
$R^+=\{\za^1,\ldots,\za^\nu\}$, where
$\za^t = s_{i_1}s_{i_2} \cdots s_{i_{t - 1}}(\alpha_{i_t})$ for
$t = 1, \ldots, \nu$.
We say that an element $\mathbf{ c} = (c_1, \ldots, c_{\nu}) \in \mathbb{
N}^{\nu}$ is of $\mathbf{ i}$-homogeneity $\mathbf{ d} =
(d_1, \ldots, d_n) \in \mathbb{ N}^n$ if
$$\sum_{t = 1}^{\nu} c_t\ \alpha^t = \sum_{k = 1}^n d_k\ \alpha_k\ .$$
Let $\mathbf{ b}(t)$ be the vector
$(0,\dots, 0, 1, 0, \dots, 0) \in \mathbb{ N}^{\nu}$ whose only
nonzero component is in the $t^{\textup{th}}$ column and is 1,
where $1 \leq t \leq \nu$.

For $\mathbf{ i}= (i_1, \ldots, i_{\nu}) \in
{\mathcal I}$ and
$\mathbf{ c} = (c_1, \ldots, c_{\nu}) \in \mathbb{ N}^{\nu}$, define
$$E_\mathbf{ i}^\mathbf{ c} = E_{i_1}^{(c_1)} \ \tilde
T_{i_1}\left(E_{i_2}^{(c_2)}\right) \ \tilde T_{i_1} \tilde
T_{i_2}\left(E_{i_3}^{(c_3)}\right) \ldots \ \tilde T_{i_1}
\tilde T_{i_2} \cdots \tilde T_{i_{(\nu -
1)}}\left(E_{i_{\nu}}^{(c_{\nu})}\right)
=\prod_{t=1}^\nu E_\mathbf{i}^{c_t \mathbf{b}(t)}.$$

Note that if $\mathbf{c}$ is of $\mathbf{ i}$-homogeneity $\mathbf{ d}$
then $\Tr(\wt(E_\mathbf{i}^\mathbf{c}))=\sum_{k=1}^n d_k$

\begin{proposition}\label{pbw}
Let $\mathbf{ i} \in {\mathcal I}$. Then
$B_\mathbf{ i} = \{E_\mathbf{ i}^\mathbf{ c} \mid \mathbf{ c} \in \mathbb{
N}^{\nu}\}$ is a $\mathbb{ Q}(v)$-basis of $\mathbf{ U}^+$.
We say that $B_\mathbf{ i}$ is a basis of PBW type. \hfill
\end{proposition}
\begin{proof}
\cite[sect. 1.8 and 1.13]{L2}
\end{proof}

We now recall Lusztig's construction of the canonical basis
of $\mathbf{ U}^+$.

\begin{theorem}\label{canonical}
Let $\mathbf{ i} \in
{\mathcal I}$ and ${\mathcal L}_\mathbf{ i}$ the $\mathbb{
Z}[v^{-1}]$-submodule of $\mathbf{ U}^+$ generated by $B_\mathbf{ i}$.
\begin{itemize}
\item[(i)] {${\mathcal L}_\mathbf{ i}$ is independent of $\mathbf{ i}$.
We denote ${\mathcal L}_\mathbf{ i}$ by ${\mathcal L}$.}
\item[(ii)] {$\pi(B_\mathbf{ i})$ is a $\mathbb{ Z}$-basis of
${\mathcal L}/v^{-1}{\mathcal L}$ independent of $\mathbf{ i}$. Here
$\pi:{\mathcal L} \rightarrow {\mathcal L}/v^{-1}{\mathcal L}$ is the
canonical projection.
We denote $\pi(B_\mathbf{ i})$ by $B$.}
\item[(iii)] { The restriction of $\pi:{\mathcal L}
\rightarrow {\mathcal L}/v^{-1} {\mathcal L}$ defines an isomorphism
of $\mathbb{ Z}$-modules $\pi^{\prime}:{\mathcal L} \cap \overline{\mathcal
L} \rightarrow {\mathcal L}/v^{-1}{\mathcal L}$ where $\overline{\mathcal
L}$
is the image of ${\mathcal L}$ under $\overline{(\ )}$. In
particular,
$\mathbf{ B} = {\pi'}^{-1}(B)$ is a $\mathbb{ Z}$-basis of ${\mathcal L}
\cap \overline {\mathcal L}$.}
\item[(iv)] {$\mathbf{ B}$ is a $\mathbb{ Z}[v^{-1}]$-basis of
${\mathcal L}$ and a $\mathbb{ Q}(v)$-basis of $\mathbf{
U}^+$. $\mathbf{ B}$ is said to be the canonical basis
of $\mathbf{ U}^+$}.
\item[(v)] {Each element of $\mathbf{ B}$ is fixed by
$\overline{(\ )}:\mathbf{ U}^+ \rightarrow \mathbf{ U}^+$.}\hfill
\end{itemize}
\end{theorem}
\begin{proof}\cite{lusztig1}
\end{proof}

\end{subsection}
\begin{subsection}{Specialization at $v=1$}

Let $\mathcal{A}=\mathbb{Q}[v]_{(v-1)}$
denote the subring of
$\mathbb{Q}(v)$ consisting of functions regular at $v=1$. Define the
$\mathcal{A}$-form $\mathbf{U}_{\mathcal{A}}$ of $\mathbf{U}^+$ to be the
$\mathcal{A}$-subalgebra of $\mathbf{U}^+$ generated by $E_i$, $1\leq i\leq
n$.
Denote by
$\mathbf{U}_1^+=\mathbf{U}_{\mathcal{A}}/(v-1)\mathbf{U}_{\mathcal{A}}$ the
specialization of $\mathbf{U}^+$ at $v=1$.
This is the positive part of the classical universal enveloping algebra
of type $\Delta$, with generators $E_1,\ldots, E_n$ and classical Serre
relations
$$\left\{\begin{array}{ll}
E_i^2 E_j - 2\, E_i E_j E_i + E_j E_i^2= 0
& \textup{if } a_{ij}=-1\\
E_i E_j - E_j E_i = 0 & \textup{if } a_{ij}\ne -1.
\end{array}\right.
$$
Note that in $\mathbf{U}_1$, $\tilde T_i(E_j)$
becomes the usual bracket $E_jE_i-E_iE_j$ and
the specialization
$B_\mathbf{i}(1)=\{E_\mathbf{i}^\mathbf{c}(1)\mid \mathbf{c }\in
\mathbb{N}^\nu\}$
of the PBW-basis is a PBW-basis of $\mathbf{U}_1$.

\end{subsection}


\begin{subsection}{Quiver modules}

Let ${\mathcal Q}=({\mathcal Q}^0,{\mathcal Q}^1)$ be a quiver whose
underlying graph
is $\Delta$, i.e. for
each edge $\{i, j\}$ of $\Delta$ we fix an orientation.
(We use the notation ${\mathcal Q}^0$ for the set of vertices of the
quiver
${\mathcal Q}$ and ${\mathcal Q}^1$ for the set of arrows.)
A vertex $i\in \mathcal{Q}^0$ is a sink (respectively a source) of
${\mathcal Q}$ if there is no
arrow $i \rightarrow j$ (respectively $i \leftarrow j$) $\in
\mathcal{Q}^1$.
An element $\mathbf{ i} = (i_1, \ldots, i_{\nu})
\in {\mathcal I}$ is {\it adapted to the quiver} ${\mathcal Q}$
if $i_1$ is a sink of ${\mathcal Q}_1 = {\mathcal Q}$
and $i_k$ is a sink of the quiver
${\mathcal Q}_k = s_{i_{k - 1}}({\mathcal Q}_{k - 1})$ obtained from
${\mathcal
Q}_{k-1}$ by reversing the orientation of all arrows
ending at
$i_{k-1}$, where $2 \leq k \leq \nu$.
It is easy to see that there is an element $\mathbf{ i}\in {\mathcal I}$
adapted to ${\mathcal Q}$ and that an element $\mathbf{ i}$ of ${\mathcal I}$
can be adapted to at most one quiver.

For the rest of this paper, let $\mathbf{ i}$ be adapted to the quiver
${\mathcal Q}$.

Let $F$ be any field. A module (or representation)
$\mathbf{ V} = (V_i, f_{ij})$ of ${\mathcal Q}$ is
a collection of $n$ finite dimensional $F$-vector spaces $V_i$, ($1\le i\le
n$)
and of $(n-1)$ $F$-linear maps $f_{ij}:V_i
\rightarrow V_j$, ($i \rightarrow j\in {\mathcal Q}^1$).
A morphism from the module $\mathbf{ V} = (V_i,
f_{ij})$ to the module $\mathbf{ V}^{\prime} = (V_i^{\prime},
f_{ij}^{\prime})$ is a collection of $F$-linear maps
$g_i:V_i \rightarrow V_i^{\prime}$, $1\le i\le n$ such that
$f_{ij}^{\prime}\circ g_i = g_j\circ f_{ij}$ for each $i
\rightarrow j \in {\mathcal Q}^1$. These modules and morphisms form
an abelian category ${\rm Mod}({\mathcal Q})$.
If $\mathbf{ V}$ is a module of ${\mathcal
Q}$, denote by $[\mathbf{ V}]$ its isomorphism class
in ${\rm Mod}({\mathcal Q})$.

The dimension of the module $\mathbf{ V} = (V_i, f_{ij})$ is the
$n$-tuple $$\dim(\mathbf{ V}) = (\dim_F(V_1), \dim_F(V_2),
\ldots,\dim_F(V_n)) \in \mathbb{ N}^n.$$
A module $\mathbf{ V}$ of ${\mathcal Q}$ is indecomposable if $\mathbf{ V}$
cannot be written as the direct sum of proper submodules.

\begin{theorem}\label{gabriel}
\begin{itemize}
\item[(i)] For all $\alpha \in R^+$, there is a unique
indecomposable module (up to isomorphism), denoted
$\mathbf{ e}_{\alpha} \in {\rm Mod}({\mathcal Q})$, such that
$\dim(\mathbf{ e}_{\alpha}) = (d_1, \ldots, d_n)$ and
$\alpha = \sum_{i = 1}^n d_i \alpha_i$; any indecomposable module
is isomorphic to $\mathbf{ e}_{\alpha}$ for a unique
$\alpha$.
This is Gabriel's theorem.
\item[(ii)] There exists an ordering $\alpha^t$, $1\leq t\leq\nu$ of the
positive roots such that ${\rm Hom}_{\mathcal Q}(\mathbf{ e}_{\alpha^t},\mathbf{
e}_{\alpha^s})=0$ if $s<t$.
\item[(iii)] There exists a bijection $\mathbf{ c} = (c_1, c_2, \ldots,
c_{\nu}) \mapsto [\mathbf{ e}(\mathbf{ c})]$ between $\mathbb{ N}^{\nu}$ and the set of
isomorphism classes of modules of ${\mathcal Q}$,
where $\mathbf{ e}(\mathbf{ c})=\oplus_{t = 1}^\nu c_t \ \mathbf{ e}_{\alpha^t}$.
In this case, $\dim(\mathbf{ e}(\mathbf{ c})) = (d_1, \ldots,
d_n)$, where $\sum_{t = 1}^{\nu} c_t\ \alpha^t =
\sum_{i = 1} d_i\ \alpha_i$, i.e. $\mathbf{ c}$ is of $\mathbf{
i}$-homogeneity $\mathbf{ d}$.
\end{itemize}
\noindent In particular, the classification of indecomposable modules is independent
of the ground field.
\end{theorem}
\begin{proof} \cite[sect. 4.12 -- 4.15]{lusztig1}
\end{proof}

Set $\left[\mathbf{ V}, \mathbf{
V}^{\prime}\right]=\dim_F {\rm Hom}_{\mathcal Q} (\mathbf{ V}, \mathbf{
V}^{\prime})$ and
$\left[\mathbf{ V}, \mathbf{ V}^{\prime}\right]^1= \dim_F {\rm Ext}_{\mathcal
Q}^1(\mathbf{ V}, \mathbf{ V}^{\prime}) $.
Note that
${\rm Hom}_{\mathcal Q}(\mathbf{ V}, \mathbf{ V}^{\prime})$ is the $F$-vector space
of morphisms $g:\mathbf{ V} \rightarrow \mathbf{ V}^{\prime}$
in ${\rm Mod}({\mathcal Q})$
and ${\rm Ext}_{\mathcal Q}^1(\mathbf{ V}, \mathbf{ V}^{\prime})$ is the $F$-vector
space of
extensions $0 \rightarrow \mathbf{ V}^{\prime} \rightarrow \mathbf{
E} \rightarrow \mathbf{ V} \rightarrow 0$ in
${\rm Mod}({\mathcal Q})$.

For
$\mathbf{ d} = (d_1, \ldots, d_n) \in
\mathbb{ N}^n$, define
$$E_\mathbf{ d} = \bigoplus_{i \rightarrow j\in{\mathcal Q}^1} {\rm
Hom}_F(F^{d_i}, F^{d_j})
\quad \textup{and} \quad G_\mathbf{ d} =\prod_{i = 1}^n GL_{d_i}(F).$$
The group $G_\mathbf{ d}$ acts
on $E_\mathbf{ d}$ by
$(g \cdot f)_{i \rightarrow j} = (g_j\
f_{ij}\ g_i^{-1})_{i\rightarrow j}.$
An element of $E_\mathbf{ d}$ can be seen as a module in ${\rm
Mod}({\mathcal Q})$ of dimension $\mathbf{ d}$. Two elements of
$E_\mathbf{ d}$ define isomorphic modules if and only
if they are in the same $G_\mathbf{ d}$-orbit. By theorem \ref{gabriel},
there
exists a bijection between the set of
$\nu$-tuples $\mathbf{ c} = (c_1, \ldots, c_{\nu})$ of
$\mathbf{ i}$-homogeneity $\mathbf{ d}$ and the set of $G_\mathbf{
d}$-orbits in $E_\mathbf{ d}$,
where $\mathbf{ c} = (c_1,
\ldots, c_{\nu})$ corresponds to the orbit ${\mathcal O}_\mathbf{
c}$ whose elements are isomorphic to $\mathbf{ e}(\mathbf{ c})$.

There is a partial order on $\mathbb{ N^\nu}$ given by
$\mathbf{ c}^{\prime} \preceq \mathbf{ c}$
if $\mathbf{ c}^{\prime}$ and $\mathbf{ c}$ have the same $\mathbf{
i}$-homogeneity and the orbit ${\mathcal O}_{\mathbf{ c}^{\prime}}$
is contained in the Zariski closure $\overline{\mathcal O}_\mathbf{
c}$ of ${\mathcal O}_\mathbf{ c}$. This is the so-called degeneration ordering.

Let $\cS$ be the set of non-split short exact sequences of modules
of ${\mathcal Q} $ and
let $\sop$ be the subset of $\cS$ consisting of all sequences for which the
first and the last module are indecomposable. Hence if $\zS\in \sop$ then
$$\zS\ :\ \suitr{\mathbf{ e}_{\za^{s}}}{\mathbf{ V}}{\mathbf{ e}_{\za^{t}}}$$
for some $s,t\in\{1,\ldots,\nu\}$ and some module $\mathbf{ V}$.
The elements of $\sop$ are called {\it elementary operations}.
For $\zS\ :\ \suitr{\mathbf{ e}_{\za^{s}}}{\mathbf{ V}}{\mathbf{ e}_{\za^{t}}} \ \in
\sop$ define
$in(\zS)=s$ and $out(\zS)=t$ and
denote by $\mathbf{ op}^{\zS}$ the vector
$(op^\zS_1,\ldots,op^\zS_\nu)\in \mathbb{ Z}^\nu$
given by
$$op^\zS_r=\left\{\begin{array}{cll}
-1&\textup{if }r=s,t\\
1&\textup{if } \mathbf{ e}_{\za^r} \textup{ is a direct summand of }
\mathbf{ V}\\
0&\textup{otherwise}.
\end{array}\right. $$
For all $\mathbf{ c}\in \mathbb{ N}^\nu$ define $\sopc{\mathbf{
c}}=\{\zS\in\sop\mid \mathbf{ c}+\mathbf{ op}^\zS\in\mathbb{ N}^\nu\}$.
Thus an elementary operation $\zS \in \sopc{\mathbf{ c}}$ allows us to go from
one orbit $\orb{c}$ to another orbit $\mathcal{O}_{\mathbf{ c}+\mathbf{ op}^\zS}$.
As we will see in Theorem \ref{Bongartz} below, elementary operations do not
only preserve the $\mathbf{ i}$-homogeneity but they are also compatible with
the partial ordering $\preceq$.

The following theorem is shown in \cite{bongartz1}.

\begin{theorem}\label{Bongartz}
Let $\mathbf{ c},\mathbf{ c}'\in \mathbb{ N}^\nu$. Then
the following four statements are
equivalent:
\begin{itemize}
\item[(i)] $ \mathbf{ c}'\preceq \mathbf{ c}$
\item[(ii)] There is a sequence of elementary operations
$\zS_1,\zS_2,\ldots,\zS_k$
such that $\zS_l\in\sop(\mathbf{ c}'+\sum_{i=1}^{l-1} \mathbf{ op}^{\zS_i})$ and
$\mathbf{ c}'+\sum_{i=1}^k \mathbf{
op}^{\zS_i} =\mathbf{ c}$
\item[(iii)] $[\mathbf{ e}_\za,\mathbf{ e}(\mathbf{ c}')]\ge [\mathbf{ e}_\za,\mathbf{ e}(\mathbf{
c})]$
for all indecomposable modules $\mathbf{ e_\za}$.
\item[(iv)] $[\mathbf{ e}(\mathbf{ c}'),\mathbf{ e}_\za]\ge [\mathbf{ e}(\mathbf{ c}),\mathbf{
e}_\za]$
for all indecomposable modules $\mathbf{ e_\za}$.
\end{itemize}
\end{theorem}

\end{subsection}


\begin{subsection}{Hall algebras}

Fix a quiver $\mathcal{Q}$ and a $\nu$-tuple $\mathbf{i}$ adapted to
$\mathcal{Q}$.
Let $\mathcal{H}_{\mathcal{Q}}=\mathcal{H}$ be the twisted Hall algebra
associated to the quiver $\mathcal{Q}$.
$\mathcal{H}$ is the free $\mathbb{Q}(v)$-module with basis
$\mathcal{B}$
the set of isomorphism classes of representations of the quiver
$\mathcal{Q}$, with multiplication
\begin{equation}
[\mathbf{e}(\mathbf{ c}')]\cdot[\mathbf{ e}(\mathbf{c}'')] =
v^{<\mathbf{e}(\mathbf{c}'),\mathbf{e}(\mathbf{c}'')>}\
\sum_{\mathbf{ c}\in \mathbb{N}^\nu} F_{\mathbf{ c}',\mathbf{
c}''}^{\mathbf{ c}}\ [\mathbf{e}(\mathbf{ c})]
\end{equation}
where
${<\mathbf{e}(\mathbf{c}'),\mathbf{e}(\mathbf{c}'')>}
=\dim_{\mathbb{F}_{v^2}}\Hom(\mathbf{e}(\mathbf{c}'),\mathbf{e}(\mathbf{c}''))
-\dim_{\mathbb{F}_{v^2}}\Ext(\mathbf{e}(\mathbf{c}'),\mathbf{e}(\mathbf{c}''))$
and
$F_{\mathbf{c}',\mathbf{c}''}^{\mathbf{c}}$
is the number
of submodules of $\mathbf{e}(\mathbf{c})$ that are isomorphic to
$\mathbf{e}(\mathbf{c}'')$ and are such that the corresponding
quotient module is isomorphic to $\mathbf{e}(\mathbf{c}')$,
as representations over $\mathbb{F}_{v^2}$.
This defines polynomials $F_{\mathbf{c}',\mathbf{c}''}^{\mathbf{c}}$
which are called Hall polynomials.

\begin{theorem}\label{Hall2}
There exists an isomorphism $\eta: \mathbf{U}^+ \to \mathcal{H}$ of
$\mathbb{Z}^n$-graded $\mathbb{Q}(v)$-algebras such that
$\eta(E_i)=[\mathbf{e}_{\za_i}]$.
It maps $E_{\mathbf{i}}^{\mathbf{c}}$ to
$v^{[\mathbf{ e}(\mathbf{ c}),\mathbf{ e}(\mathbf{ c})]
-\Tr(\wt(E_{\mathbf{i}}^{\mathbf{c}}))}
[\mathbf{ e}(\mathbf{ c})]$.
\end{theorem}
\begin{proof}
\cite{ringel4}\cite{green1}
\end{proof}

We shall need the following corollary, see \cite{caldero1} :
\begin{cor}\label{filtration}
Up to a power of $v$, the
$E_{\mathbf{i}}^{\mathbf{c'}}$-coefficient of $\prod_{i=1}^k
E_{\mathbf{i}}^{\mathbf{c}_i}$ is
$F_{\mathbf{c}_1,\ldots,\mathbf{c}_k}^{\mathbf{c}'}(v^2)$, where the polynomial
$F_{\mathbf{c}_1,\ldots,
\mathbf{c}_k}^{\mathbf{c'}}(q)$
denotes the number of filtrations of $\mathbf{e}
(\mathbf{c'})$ with successive quotients isomorphic to
$\mathbf{e}(\mathbf{c}_1)$, $\ldots$, $\mathbf{e}(\mathbf{c}_k)$ over
$\mathbb{F}_{q}$.

\end{cor}
\end{subsection}

\begin{subsection}{Coalgebra structure}
Define an algebra structure on $\mathbf{U}^+\otimes_{\mathbb{Q}(v)}
\mathbf{U}^+$ by $$(x\otimes x').(y\otimes y')=v^{(\wt(x'),
\wt(y))}(xy)\otimes(x'y')$$ on homogeneous elements. Let
$\delta:\mathbf{U}^+\rightarrow \mathbf{U}^+\otimes_{\mathbb{Q}(v)}
\mathbf{U}^+$ be the $\mathbb{Q}(v)$-algebra map given by
$\delta(E_i)=E_i\otimes 1+1\otimes E_i$. Now, there exists a unique
$\mathbb{Q}(v)$-bilinear form $(\,,\,)$ on $\mathbf{U}^+$ such that
$(E_i,E_j)=\delta_{i,j}(1-v^{-2})^{-1}$ and $(x,yy')=(\delta(x),y\otimes
y')$, where $(\,,\,)$ is extended to $\mathbf{U}^+\otimes \mathbf{U}^+$ by
the rule $(x\otimes x',y\otimes y')=(x,y).(x',y')$.

The Hall algebra approach of the quantum algebra $\mathbf{U}^+$ gives a nice
description of its coalgebra structure, \cite{green1}.
For every
representation $\mathbf{e}(\mathbf{c})$ of the quiver $\mathcal{Q}$, let
$a_{\mathbf{c}}(v^2)$ be the number of automorphisms of
$\mathbf{e}(\mathbf{c})$ as a
representation over $\mathbb{F}_{v^2}$. Recall \cite{Ringel3} that
\begin{equation}\label{automorphism} a_{\mathbf{c}}(v^2)=
v^{2\sum_{s<t}c_sc_t[\mathbf{e}_{\alpha^s},\mathbf{e}_{\alpha^t}]}
\prod_{s=1}^{\nu}\mid GL_{c_s}(\mathbb{F}_{v^2}))\mid
\end{equation}
We have \cite{reineke1}:
\begin{equation}\label{orthogonal}(E_{\mathbf{i}}^{\mathbf{c}},E_{\mathbf{i}}^{\mathbf{c'}})=
\delta_{\mathbf{c},\mathbf{c'}}\,
v^{2[\mathbf{e}(\mathbf{c}),\mathbf{e}(\mathbf{c})]}\,
a_{\mathbf{c}}(v^2)^{-1}.\end{equation}
This proves that the basis $B_{\mathbf{i}}$ is orthogonal. We can define the
dual basis for this pairing
$B_{\mathbf{i}}^*=\{E_{\mathbf{i}}^{\mathbf{c}*},
\mathbf{c}\in\mathbb{N}^\nu\}$, with
\begin{equation}\label{dualbasis}E_{\mathbf{i}}^{\mathbf{c}*}=
v^{-2[\mathbf{e}(\mathbf{c}),\mathbf{e}(\mathbf{c})]}
\,a_{\mathbf{c}}(v^2\,)E_{\mathbf{i}}^{\mathbf{c}}.
\end{equation}
This implies that the algebra $\mathbf{U}^+$ is generated by
$E_{\mathbf{i}}^{\mathbf{b}(s)*}$, $1\leq s\leq\nu$, with straigh\-tening
relations :
\begin{equation}\label{straightening}
E_{\mathbf{i}}^{\mathbf{b}(t)*}
E_{\mathbf{i}}^{\mathbf{b}(s)*}=
\sum_{\mathbf{b}(t)+\mathbf{b}(s)\preceq\mathbf{c}}
G_{\mathbf{b}(t),\mathbf{b}(s)}^{\mathbf{c}}
E_{\mathbf{i}}^{\mathbf{c}*}.
\end{equation}
In order to understand those straightening relations, we have to define
another polynomial. Let
$0\rightarrow\mathbf{e}(\mathbf{c}')\rightarrow
\mathbf{e}(\mathbf{c})\rightarrow
\mathbf{e}(\mathbf{c}'')\rightarrow 0$
be a non-split exact sequence of representations of
$\mathcal{Q}$. Then, the set of points
$\mathcal{C}_{\mathbf{c}'',\mathbf{c}'}^{\mathbf{c}}$ in the
$\mathbb{F}_{v^2}$-space Ext$^1(\mathbf{e}(\mathbf{c}''),
\mathbf{e}(\mathbf{c}'))$ corresponding to
this extension is a cone (minus the zero point). Set
$$E_{\mathbf{c}'',\mathbf{c}'}^{\mathbf{c}}(v^2):=\ \mid\mathcal{C}_{\mathbf{c}
'',\mathbf{c}'}^{\mathbf{c}}\mid.$$
It is known that
$E_{\mathbf{c}'',\mathbf{c}'}^{\mathbf{c}}$ is a polynomial.

\begin{proposition}\label{riedtmann} Suppose $1\leq s<t\leq\nu$. Then, for
all $\mathbf{c}$ in $\mathbb{N}^\nu$ and up to a power of $v$, we have
$$
G_{\mathbf{b}(t),\mathbf{b}(s)}^{\mathbf{c}} =
\left\{\begin{array}{ll}
1 & \textup{if } \mathbf{c}=\mathbf{b}(s)+\mathbf{b}(t),\\
E_{\mathbf{b}(t),\mathbf{b}(s)}^{\mathbf{c}}(v^2) & \textup{if }
\mathbf{c}=\mathbf{b}(s)+\mathbf{b}(t)+\mathbf{op}^\zS \textup{ for some }
\zS \in \sop ,\\
0 &\textup{otherwise.}
\end{array}\right.$$
\end{proposition}

\begin{proof}
This is a direct consequence of (\ref{dualbasis}) and the Riedtmann formula
which states
\begin{equation}
F_{\mathbf{b}(t),\mathbf{b}(s)}^{\mathbf{c}}(v^{2})\,
a_{\mathbf{b}(t)}(v^{2})\,a_{{\mathbf{b}(s)}}(v^{2}) \,
a_{\mathbf{c}}(v^{2})^{-1}\,=\,
E_{\mathbf{b}(t),\mathbf{b}(s)}^{\mathbf{c}}(v^{2}).
\end{equation}
\end{proof}

\begin{remark} As the set
$\mathcal{C}_{\mathbf{c}'',\mathbf{c}'}^{\mathbf{c}}\cup\{0\}$ is a cone, we
see that the polynomial $E_{\mathbf{c}'',\mathbf{c}'}^{\mathbf{c}}(v^2)$ can
be factorized by $v^2-1$. So, the proposition shows that the
$\mathcal{A}$-space generated by the $E_{\mathbf{i}}^{\mathbf{c}*}$ is an
algebra, and its specialization at $v=1$ is a (commutative) polynomial
algebra. Indeed, it is known that this specialization is the algebra of
regular functions on the maximal unipotent group.
\end{remark}

Now, let $\sigma$ be the $\mathbb{Q}$-antiautomorphism of $\mathbf{U}^+$
such that $\sigma(E_i)=E_i$, $1\leq i\leq n$, $\sigma(v)=v^{-1}$.
\begin{lemma}\label{adjoint}
Let $x$, $y$ be homogeneous elements in $\mathbf{U}^+$, and let
$t\in\{1,\ldots,\nu\}$. We have
$$ \begin{array}{crcl}
\textup{(i)}&(\overline x,y)&=&(-v)^{\Tr(\wt(x))}\,v^{{1\over 2}(\wt(x),
\wt(x))_Q}\,\overline{(x,\sigma(y))}.\\
\textup{(ii)}&\sigma(E_{\mathbf{i}}^{\mathbf{b}(t)})&=&(-1)^
{\Tr(\alpha^t)-1}\,E_{\mathbf{i}}^{\mathbf{b}(t)}
,\textup{ up to a power of $v$}.
\end{array}$$
\end{lemma}
\begin{proof}
(i) is a direct consequence of~\cite[4.3]{reineke1} and (iii)
is proved in~\cite[prop. 2.1]{caldero1}.
\end{proof}

\end{subsection}


\begin{subsection}{Local intersection cohomology of orbit closures}
In this subsection, let $\mathbb{F}$ be  an algebraic closure
of a finite field $\mathbb{F}_q$ with $q=p^e$ elements, where $p$ is a prime number,
and let $\mathbf{ d} = (d_1, \ldots, d_n)
\in \mathbb{ N}^n$. We will write the dimension dim$(\mathcal{O}_{\mathbf{c}})$ of the orbit $\mathcal{O}_{\mathbf{c}}$ by $d(\mathbf{c})$.

The results of this subsection have been proved in \cite[chap. 9 --
10]{lusztig1}.

\begin{proposition}\label{omega}
Let $\mathbf{ c}\in\mathbb{ N}^\nu$ be of $\mathbf{ i}$-homogeneity $\mathbf{
d}$. Then for each $\mathbf{ c}'\preceq \mathbf{ c}$, there exists $\omega_\mathbf{
c^{\prime}}^\mathbf{ c} \in \mathbb{ Z}[v, v^{-1}]$ such that
$$\overline{E_\mathbf{ i}^\mathbf{ c}} =
\sum_{\mathbf{ c}^{\prime}\preceq\mathbf{ c}}
\omega_\mathbf{ c^{\prime}}^\mathbf{ c}\ E_\mathbf{ i}^\mathbf{
c^{\prime}}.$$
Moreover $\omega_\mathbf{ c}^\mathbf{ c} = 1$ and
for all $\mathbf{ c}'\preceq \mathbf{ c}$, $\Omega_\mathbf{ c^{\prime}}^\mathbf{ c} \soit
v^{d(\mathbf{ c}) - d(\mathbf{ c^{\prime}})} \omega_\mathbf{
c^{\prime}}^\mathbf{ c}$ is an element of $\mathbb{ Z}[v^2,v^{-2}]$.
\end{proposition}

\begin{theorem}\label{zeta}
{ Let $\mathbf{ c}\in \mathbb{ N}^{\nu}$ be of $\mathbf{ i}$-homogeneity
$\mathbf{ d}$ and let ${\mathcal E}^\mathbf{ c}\in\mathbf{ B}$ be the unique canonical
basis element
such that $\pi({\mathcal E}^\mathbf{ c}) = \pi(E_\mathbf{
i}^\mathbf{ c})$.} Then
\begin{itemize}
\item[(i)] { ${\mathcal E}^\mathbf{ c} =
\sum_\mathbf{ c^{\prime}} \zeta_\mathbf{ c^{\prime}}^\mathbf{ c}\ E_\mathbf{
i}^\mathbf{ c^{\prime}}$, where
$\mathbf{ c}^{\prime}$ runs over the set of elements of $\mathbb{
N}^{\nu}$ of $\mathbf{ i}$-homogeneity $\mathbf{ d}$, $\zeta_\mathbf{
c}^\mathbf{ c} = 1$ and $\zeta_\mathbf{ c^{\prime}}^\mathbf{ c} \in
v^{-1}\mathbb{ Z}[v^{-1}]$ for $\mathbf{ c}^{\prime} \ne \mathbf{
c}$.}
\item[(ii)] {If $\mathbf{ c}^{\prime} \not \preceq \mathbf{ c}$, then
$\zeta_\mathbf{ c^{\prime}}^\mathbf{ c} = 0$.}
\item[(iii)] {If $\overline{(\ )}$ is the $\mathbb{ Z}$-linear involution
of $\mathbb{ Z}[v , v^{-1}]$ sending $v$ to
$v^{-1}$, then $$\zeta_\mathbf{ c^{\prime}}^\mathbf{ c} =
\sum_{\mathbf{ c}^{\prime\prime}\atop \mathbf{ c}^{\prime} \preceq
\mathbf{ c}^{\prime\prime} \preceq \mathbf{ c}} w_\mathbf{
c^{\prime}}^\mathbf{ c^{\prime\prime}}\ \overline{\zeta_\mathbf{
c^{\prime\prime}}^\mathbf{ c}}.$$}
\item[(iv)] {If $\mathbf{ c}^{\prime} \preceq \mathbf{ c}$, $f$ is a
$\mathbb{F}_q$-rational point of the orbit ${\mathcal O}_\mathbf{
c^{\prime}}$ in $E_\mathbf{ d}$
and ${\mathcal H}_f^a$ is the stalk at $f$ of the $a$-th
cohomology sheaf of the
intersection cohomology complex of the Zariski closure $\overline{\mathcal
O}_\mathbf{ c}$ of ${\mathcal O}_\mathbf{ c}$ with coefficients in
$\overline{\mathbb{ Q}_{\ell}}$ (extended by zero on the complement of
that closure), where $\ell$ is a prime number $\ne p$,
and with the $\mathbb{F}_q$-structure such that the Frobenius map acts
as identity on the stalks of its 0-th cohomology sheaf
at the rational points of the orbit
${\mathcal O }_\mathbf{ c}$,
then
$${\mathcal H}_f^{2a + 1} = 0\quad \hbox{ for all } a \quad
\hbox{ and } \quad v^{d(\mathbf{ c}) - d(\mathbf{ c^{\prime}})}
\zeta_\mathbf{ c^{\prime}}^\mathbf{ c} = \sum_a \dim({\mathcal
H}_f^{2a})\ v^{2a}.$$
In particular, $v^{d(\mathbf{ c}) - d(\mathbf{
c}^{\prime})} \zeta_\mathbf{ c^{\prime}}^\mathbf{ c}$
is a polynomial in $v^2$ with coefficients in $\mathbb{ N}$. }
\end{itemize}
\end{theorem}

\begin{definition}\label{rsdef}
We say that the orbit closure $\overline {\mathcal O}_{\mathbf{ c}}$ is
{\it rationally smooth at} ${\mathcal O}_{\mathbf{ c}'}\subset \overline
{\mathcal O}_{\mathbf{ c}}$ if
for all $\mathbf{ c}''$ such that $\mathbf{ c}'\preceq \mathbf{
c}''\preceq \mathbf{ c}$ we have
$\sum_a\dim({\mathcal H}_f^{2a})v^{2a}=1$ for a $\mathbb{F}_q$-rational point $f\in
{\mathcal O}_{\mathbf{ c}''}$, i.e. if $\zeta_{\mathbf{ c}''}^{\mathbf{ c}}=v^{d(\mathbf{
c''})-d(\mathbf{ c})}$ for all $\mathbf{ c}''$ such that $\mathbf{ c}'\preceq \mathbf{
c}''\preceq \mathbf{ c}$.\\
The orbit closure $\overline {\mathcal O}_{\mathbf{ c}}$ is
{\it rationally smooth } if it is rationally smooth at each ${\mathcal
O}_{\mathbf{ c}'}\subset\overline
{\mathcal O}_{\mathbf{ c}}$, i.e. if $\zeta_{\mathbf{ c}'}^\mathbf{ c}=v^{{d(\mathbf{
c'})-d(\mathbf{ c})}}$ for all $\mathbf{ c}'\preceq \mathbf{ c} $.
\end{definition}

\begin{remark}\label{bormac}
  By a result of \cite{bormac}, $\overline {\mathcal O}_{\mathbf{ c}}$ is
  rationally smooth  iff for each point $x\in\overline {\mathcal
  O}_{\mathbf{ c}}$ we have
  $$H^i(\overline {\mathcal O}_{\mathbf{ c}},\overline {\mathcal
  O}_{\mathbf{ c}}\setminus \{x\};\mathbb{Q}) = \left\{
  \begin {array}{ll}
  \mathbb{Q} & \textup{if } i= 2\dim \overline {\mathcal O}_{\mathbf{ c}}\\
  0          & \textup{otherwise}
  \end{array}\right.
  $$
  where $H^i$ denotes ordinary cohomology.
\end{remark}

\end{subsection}
\end{section}


\begin{section} {Algebraic approach}
This section gives a generalization of results in~\cite{bedschi1}. We first
prove some formulas on the coefficients $\Omega_{\mathbf{c'}}^{\mathbf{c}}$.
\begin{subsection}{A formula for $\Omega_{\mathbf{c'}}^{\mathbf{c}}$.}
We fix a quiver $\mathcal{Q}$ and a $\nu$-tuple $\mathbf{i}$ adapted to
$\mathcal{Q}$ as in the previous section.
For all $\mathbf{c}$ in $\mathbb{N}^\nu$, set
$\mathbf{c}^t:=c_t\mathbf{b}(t)$, $1\leq t\leq \nu$.
An explicit formula for $\Omega_{\mathbf{c'}}^{\mathbf{c}}$ can be
calculated :

\begin{proposition}\label{formula} For every $\mathbf{c}$, $\mathbf{c'}$ in
$\mathbb{N}^\nu$, we have :
\begin{equation} \Omega_{\mathbf{c'}}^{\mathbf{c}}=
F_{\mathbf{c'}^\nu,\ldots,\mathbf{c'}^1}^{\mathbf{c}}(v^{-2})
\left(\prod_{t=1}^\nu a_{\mathbf{{c}'}^t}(v^{-2})\right)
a_{\mathbf{c}}(v^{-2})^{-1}
\end{equation}
up to a power of $v$.
\end{proposition}

\begin{proof}
In the following proof, equalities are up to a power of $v$.
We first write
$ \Omega_{\mathbf{c'}}^{\mathbf{c}}=(\overline{E_{\mathbf{i}}^{\mathbf{c}}},
E_{\mathbf{i}}^{\mathbf{c'}*})$.
Hence, by Lemma~\ref{adjoint} (i), we have
$\overline{\Omega_{\mathbf{c'}}^{\mathbf{c}}}=(-1)^
{\Tr(\wt(E_{\mathbf{i}}^{\mathbf{c}'})
)}(E_{\mathbf{i}}^{\mathbf{c}},
\sigma(E_{\mathbf{i}}^{\mathbf{c'}*}))$.
So, by formula (\ref{orthogonal}), it remains to find the coefficient of
$E_{\mathbf{i}}^{\mathbf{c}}$ in the decomposition of
$\sigma(E_{\mathbf{i}}^{\mathbf{c'}*})$ in the basis $B_{\mathbf{i}}$.

It is easily seen that
$E_{\mathbf{i}}^{\mathbf{c'}*}=\prod_{t=1}^\nu
E_{\mathbf{i}}^{c'_t\mathbf{b}(t)*}$, up to a power of $v$.
Then, using equation (\ref{dualbasis}), Lemma~\ref{adjoint} (ii) and the
fact that
$\overline{a_{c_t\mathbf{b}(t)}(v^2)} \, {a_{c_t\mathbf{b}(t)}(v^2)}^{-1}
=(-1)^{c_t}$ up to a power of $v$, we get
\begin{equation}\label{dualequation}
\overline{\Omega_{\mathbf{c'}}^{\mathbf{c}}}
=(E_{\mathbf{i}}^{\mathbf{c}},\prod_{t=\nu}^1
E_{\mathbf{i}}^{c'_t\mathbf{b}(t)*}) \textup{ up to a power of $v$}.
\end{equation}

Finally, using equation~\ref{dualbasis}, Corollary~\ref{filtration} (i) and
equation~\ref{orthogonal}, we obtain :
$$\overline{\Omega_{\mathbf{c'}}^{\mathbf{c}}}=
F_{\mathbf{c'}^\nu,\ldots,\mathbf{c'}^1}^{\mathbf{c}}(v^2)
\left(\prod_{t=1}^\nu
a_{\mathbf{c}{'}^t}(v^{2})\right)
a_{\mathbf{c}}(v^2)^{-1}$$
and the proposition follows.
\end{proof}

There is a nice particular case of this formula when
$\mathbf{e}(\mathbf{c'})$ is the sum of two indecomposable modules.

\begin{cor}\label{elementary} Suppose
$\mathbf{c'}=\mathbf{b}(s)+\mathbf{b}(t)$, $s<t$.
Then,
up to a power of $v$,
\begin{equation} \Omega_{\mathbf{c'}}^{\mathbf{c}}=E_{\mathbf{b}(t),
\mathbf{b}(s)}^{\mathbf{c}} (v^{-2}).
\end{equation}
\end{cor}

\begin{proof}
This is a direct application of the Riedtmann formula,~\cite{riedtmann1}.
\end{proof}
\end{subsection}

\begin{subsection}{Derivative of $\Omega_{\mathbf{c'}}^{\mathbf{c}}$ at
$v=1$.}
Set
$$e_{\mathbf{c}'',\mathbf{c}'}^{\mathbf{c}}(v^2)
:={E_{\mathbf{c}'',\mathbf{c}'}^{\mathbf{c}}(v^2)\over v^2-1},$$
which is a polynomial.
Let $\sopr$ be the set of {\it regular elementary operations} defined by
$$\sopr=\{\zS\ :\ \suitr{\mathbf{e}_{\za^s}}{{\mathbf e}(\mathbf{c})}{\mathbf {e}_{\za^t}},\
e_{ {\mathbf b}(t),\mathbf{ b}(s)}^{\mathbf{ c}}(1)\not=0\}.$$
With the notation above, we set $e^{\zS}:=e_{ \mathbf{b}(t),\mathbf{b}(s)}
^{\mathbf{c}}(1)$.\par
By~\cite{norenberg1}, we know that all elementary operations are regular in
case A and D. This is not true in case E, see remark below.
For $\mathbf{c}$ in $\mathbb{N}^\nu$, define $Opr(\mathbf{c}):=Opr\cap
Op(\mathbf{c})$.
We can now generalize \cite[Theorem 5.4]{bedschi1}.

\begin{theorem} \label{main}
Let $\mathbf{c}$, $\mathbf{c'}$ in $\mathbb{N}^\nu$ with same
$\mathbf{i}$-homogeneity, then
\begin{itemize}
\item[(i)]
$\Omega_{\mathbf{c}}^{\mathbf{c}}=1$ and
$\Omega_{\mathbf{c'}}^{\mathbf{c}}(1)=0$ if ${\mathbf{c'} \prec \mathbf{c}}.$
\item[(ii)]
$\left.{{\rm d}\over {\rm dv}}\zO_{\mathbf{c}'}^{\mathbf{c}}\right|_{v=1}\not=0$ if and only
if there exists $\zS$ in $\sopr(\mathbf{c})$ such that
$\mathbf{c}=\mathbf{c}'+ \mathbf{op}^{\zS}$.
\item[(iii)] If the condition in (ii) is verified, then
$\left.{{\rm d}\over {\rm dv}}\zO_{\mathbf{c}'}^{\mathbf{c}}\right|_{v=1}
=-2\,c'_sc'_t\,e^{\zS}$,
where $s=in(\zS)$ and $t=out(\zS)$.
\end{itemize}
\end{theorem}
\begin{proof}
The first formula of (i) is given in Proposition \ref{omega}.
For the second formula, just remark that the map $\overline{(\,)}$
specializes on identity at $v=1$.
\par
Let's prove (ii) and (iii).
By equation (\ref{dualequation})
$\overline{\zO_{\mathbf{c}'}^\mathbf{c}}
= \left(E_\mathbf{i}^\mathbf{c}, \prod_{t=\nu}^1
E_\mathbf{i}^{c'_t \mathbf{b}(t)*}\right)$,
up to a power of $v$. Say
$\sum_{\mathbf{c}''}
h^{\mathbf{c}''}_{\mathbf{c}'}E_\mathbf{i}^{\mathbf{c}''*}
= \prod_{t=\nu}^1 E_\mathbf{i}^{c'_t \mathbf{b}(t)*} $
is the expansion of $ \prod_{t=\nu}^1 E_\mathbf{i}^{c'_t \mathbf{b}(t)*}$
in the dual PBW-basis,
thus $\overline{\zO_{\mathbf{c}'}^\mathbf{c}}=
h^{\mathbf{c}}_{\mathbf{c}'}$ up to a power of $v$.
Then, since
${\zO_{\mathbf{c}'}^\mathbf{c}}(1)=0$, we have
${\textup{d}\over \textup{dv}}
\overline{\zO_{\mathbf{c}'}^\mathbf{c}}|_{v=1} =
{\textup{d}\over \textup{dv}}
h^{\mathbf{c}}_{\mathbf{c}'}|_{v=1}$,
which is non-zero iff the multiplicity of $(v-1)$ in $h_{\mathbf{c}'}^\mathbf{c}$ is
one.
The above expansion is obtained by performing a certain number of
straightening relations (\ref {straightening}). By proposition
\ref{riedtmann} we have
\begin{enumerate}
\item $h_{\mathbf{c}'}^{\mathbf{c}'}=1$, up to a power of $v$.
\item If $\mathbf{c}=\mathbf{c}'+\mathbf{op}^\zS$ for some $\zS\in
\sop$ with $s=in(\zS), \ t=out(\zS)$,
then
$$h_{\mathbf{c}'}^{\mathbf{c}} =
E_{\mathbf{b}(t),\mathbf{b}(s)}^\mathbf{c} (v^2)
= (v^2-1) \, e_{\mathbf{b}(t),\mathbf{b}(s)}^\mathbf{c}(v^2) $$
up to a power of $v$. Then the multiplicity of $(v-1) $ in
$h_{\mathbf{c}'}^\mathbf{c} $
is at least one and is equal to one iff $\zS$ is regular.
\item All other
$h_{\mathbf{c}'}^\mathbf{c''} $
have $(v-1)$ multiplicity at least 2.
\end{enumerate}
This proves (ii). Now let $\mathbf{c}=\mathbf{c}'+\mathbf{op}^\zS$ with
$\zS\in \sopr$ and $s=in(\zS)$, $t=out(\zS)$. Then the only way
$E_\mathbf{i}^{\mathbf{c}*}$ terms appear in our expansion is in the
straightening relations involving both $E_\mathbf{i}^{\mathbf{b}(t)*}$
and $E_\mathbf{i}^{\mathbf{b}(s)*}$.
The number of these straightening relations is $c'_s c'_t$.
Thus, by Proposition \ref{riedtmann}, $h_{\mathbf{c}'}^\mathbf{c}$ is a sum of $c'_s c'_t$ terms equal to $(v^2-1)\,
e_{\mathbf{b}(t),\mathbf{b}(s)}^{\mathbf{c}}(v^2)$ up to a power of $v$. Hence,
${\textup{d}\over \textup{dv}}h_{\mathbf{c}'}^\mathbf{c}|_{v=1}= 2\,c'_s
c'_t\,e^\zS$.
The theorem follows since
${\textup{d}\over \textup{dv}}\zO_{\mathbf{c}'}^\mathbf{c}|_{v=1}
=-{\textup{d}\over
\textup{dv}}\overline{\zO_{\mathbf{c}'}^\mathbf{c}}|_{v=1}$.

\end{proof}
\begin{remark}
The formula in (iii) slightly differs from the one in~\cite{bedschi1}
because of the factor $e^{\zS}$. This factor can be obtained from the
polynomials calculated in~\cite{norenberg1}. In type A, it can only be 1, so
this agrees with the formula in~\cite{bedschi1}. In type D this constant can
be equal to 1 or -1. In type E, it can be equal to 1, -1 and 0.
\end{remark}

\end{subsection}

\begin{subsection}{Rational Smoothness}

In this section, we will characterize which orbit closures
$\overline{\mathcal{ O}_{\mathbf c}}$ are rationally smooth. As a
consequence, we will show that if $\overline{\mathcal{ O}_{\mathbf c}}$
is rationally smooth, then $\overline{\mathcal{ O}_{\mathbf c}}$ is smooth.

Let $u=v^2$.
The next proposition has been shown in \cite[Cor 6.4]{bedschi1}. We sketch a proof of it for completeness. \par
Suppose that $\overline{\mathcal{O}_\mathbf{c}}$ is rationally smooth, then
for all $\mathbf{c}'\preceq\mathbf{c}$ we have $\zeta_{\mathbf{c}'}^{\mathbf{c}}=v^{d(\mathbf{c}')-d(\mathbf{c})}$. Fix now $\mathbf{c}'\preceq\mathbf{c}$. We obtain by Theorem \ref{zeta} (iii) that $$\sum_{\mathbf{c}'\preceq\mathbf{c}''\preceq\mathbf{c}}\Omega_{\mathbf{c}'}^{\mathbf{c}''}u^{d(\mathbf{c})-d(\mathbf{c}'')}=1.$$
By taking the derivative to $u$ evaluated at $u=1$, we get by Theorem \ref{main} (i) :
\begin{proposition}\label{corollary 6.4}
$\overline{\mathcal{O}_\mathbf{c}}$ is rationally smooth then
$$
- \sum_{\mathbf{ c}'' \atop \mathbf{ c}' \prec \mathbf{ c}'' \preceq
\mathbf{ c}}
\left.{\textup{d} \Omega_{\mathbf{
c}'}^{\mathbf{c}''}\over \textup{du}}\right|_{u = 1}
= d(\mathbf{ c}) - d(\mathbf{ c}')
$$
for all $\mathbf{c}'\preceq \mathbf{c}$.
\end{proposition}
For a subset $J$ of the set $\mathcal{Q}^1$ of arrows of the
quiver $\mathcal{Q}$, we define
$$
E_\mathbf{ d}(J) =
\{(f_{ij})_{i\rightarrow j} \in E_\mathbf{ d} \mid f_{ij} =
0\quad \textup{if $\{i, j\} \in \mathcal{Q}^1 \setminus J$}\}.
$$
The following properties are easily proved and left to the
reader. Let $J, J^{\prime}$ be two subsets of the
set $\mathcal{Q}^1$ of edges of $\mathcal{Q}$.
\begin{enumerate}
\item $E_\mathbf{ d}(J)$ is a linear subspace of $E_\mathbf{
d}$ of dimension $\dim(E_\mathbf{ d}(J)) = \sum_{\{i,j\} \in J}
d_i\, d_j$. In particular, $E_\mathbf{ d}(J)$ is a smooth
variety.

\item $E_\mathbf{ d}(J)$ is a $G_\mathbf{ d}$-stable closed
subset of $E_\mathbf{ d}$ and it is a finite union of $G_\mathbf{
d}$-orbits. As a consequence and because the field $\mathbb{F}$ is
algebraically closed, we get that there is a unique open dense
$G_\mathbf{ d}$-orbit in $E_\mathbf{ d}(J)$. We will denote this
orbit by $\mathcal{ O}(J)$.

\item $E_\mathbf{ d}(J) \cap E_\mathbf{ d}(J') =
E_\mathbf{ d}(J \cap J')$ and $E_\mathbf{ d}(J)
+ E_\mathbf{ d}(J') = E_\mathbf{ d}(J \cup
J')$.

\item $E_\mathbf{ d}(J) \subseteq E_\mathbf{ d}(J')$
if $J \subseteq J'$.

\item $E_\mathbf{ d}(\emptyset) = \{0\}$ and $E_\mathbf{
d}(\mathcal{Q}^1) = E_\mathbf{ d}$.
\end{enumerate}

Because of (3) and (5) we see that for each $G_\mathbf{
d}$-orbit $\mathcal{ O}_\mathbf{ c}$ in $E_\mathbf{ d}$, there is a
unique smallest subset $J(\mathbf{ c})$ of $\mathcal{Q}^1$ for which
$\mathcal{ O}_\mathbf{ c} \subseteq E_\mathbf{ d}(J(\mathbf{ c}))$.
In fact
$$
J(\mathbf{ c}) = \bigcap_{\mathcal{ O}_\mathbf{ c} \subseteq E_\mathbf{
d}(J)} J.
$$

Because $E_\mathbf{ d}(J)$ is closed, then $\mathcal{ O}_\mathbf{ c}
\subseteq E_\mathbf{ d}(J)$ if and only if $\overline{\mathcal{
O}_\mathbf{ c}} \subseteq E_\mathbf{ d}(J)$.
As a consequence, if
$\mathbf{ c}'' \preceq \mathbf{ c}$, then we get
easily that $J(\mathbf{ c}'') \subseteq J(\mathbf{
c})$.
Note also that we don't necessarily have $\mathcal{
O}_\mathbf{ c} = \mathcal{ O}(J(\mathbf{ c}))$.
In fact, we will prove
in theorem~\ref{smoothness} that we get this equality precisely when
$\overline{\mathcal{ O}_\mathbf{ c}}$ is rationally smooth.

Let $\mathbf{ c}^{\min} = (c_1^{\min}, c_2^{\min}, \dots,
c_{\nu}^{\min}) \in \mathbb{ N}^{\nu}$ be defined by
$$
c_k^{\min} = \left\{ \begin{array}{ll}
d_i &\textup{if $\za^k =
\alpha_i$ for some $i$, $1 \leq i \leq n$}\\
0 &\textup{otherwise}.\end{array}\right.
$$
$\mathcal{ O}_{\mathbf{ c}^{\min}}$ is the minimal $G_\mathbf{
d}$-orbit of $E_\mathbf{ d}$ and it consists of only the identity
element of $E_\mathbf{ d}$.

\begin{theorem}\label{smoothness}
\begin{enumerate}
\item We have the equality
$$ - \sum_{\mathbf{ c}''\atop \mathbf{ c}^{\min} \prec
\mathbf{ c}'' \preceq \mathbf{ c}}
\left. {\textup{d}
\Omega_{\mathbf{ c}^{\min}}^{\mathbf{c}''} \over \textup{du}}\right|_{u = 1}
= \dim(E_\mathbf{ d}(J(c))). $$
\item
$\overline{\mathcal{ O}_\mathbf{c}}$ is rationally smooth
if and only if
$\overline{\mathcal{ O}_\mathbf{ c}} = E_\mathbf{ d}(J(\mathbf{
c}))$.
\end{enumerate}
\end{theorem}

\begin{proof}
By theorem \ref{main},
$$
- \sum_{\mathbf{ c}'' \atop \mathbf{ c}^{\textup{min}} \prec \mathbf{ c}'' \preceq\mathbf{ c}}
\left.{\textup{d} \Omega_{\mathbf{c}^{\textup{min}}}^{\mathbf{c}''}\over
\textup{du}}\right|_{u = 1}
= \sum_{\zS\in \sopr(\mathbf{c}^{\textup{min}})\atop
\mathbf{ c}^{\textup{min}}+\mathbf{op}^\zS \preceq \mathbf{ c}}
c_{in(\zS)}^{\textup{min}}\,c_{out(\zS)}^{\textup{min}} \,e^\zS.$$
Now $\zS\in \sopr(\mathbf{c}^{\textup{min}})$ iff $\zS\in \sopr$ and
$c_{in(\zS)}^{\textup{min}}\ne 0 \ne c_{out(\zS)}^{\textup{min}}$,
i.e. $\za^{in(\zS)}=\za_j$ and $\za^{out(\zS)}=\za_i$ are simple roots,
$c_{in(\zS)}^{\textup{min}}=d_j$, $c_{out(\zS)}^{\textup{min}}=d_i$
and there is an arrow $i\to j \in \mathcal{Q}^1$.
Moreover $\zS\ :\
\suitr{\mathbf{e}_{\za_j}}{\mathbf{e}_{\za_i+\za_j}}{\mathbf{e}_{\za_i}}$
and $e^\zS=1$.
Note also that in this situation $\mathbf{c}^{\textup{min}}+\mathbf{op}^\zS
\preceq \mathbf{c}$ iff $i\to j \in J(\mathbf{c})$.
Thus
$$
- \sum_{\mathbf{ c}'' \atop \mathbf{ c}^{\textup{min}} \prec \mathbf{ c}'' \preceq
\mathbf{ c}}
\left.{\textup{d} \Omega_{\mathbf{c}^{\textup{min}}}^{\mathbf{c}''}\over
\textup{du}}\right|_{u = 1}
= \sum_{i\to j \in J(\mathbf{c})} d_id_j
= \dim (E_\mathbf{d}(J(\mathbf{c}))),
$$
this proves (1).

(2) $E_\mathbf{ d}(J(\mathbf{ c}))$ is smooth, hence rationally smooth.

Conversely if $\overline {\mathcal{ O}_\mathbf{ c}}$ is rationally
smooth, then applying proposition \ref{corollary 6.4} for $\mathbf{
c}' = \mathbf{ c}^{\min}$ and using part (1), we get
$\dim(E_\mathbf{d}(J(\mathbf{c}))) = d(\mathbf{ c}).$
Hence $\mathcal{O}_\mathbf{c} $ is the unique dense orbit in
$E_\mathbf{d}(J(\mathbf{c}))$, thus
$\overline {\mathcal{ O}_\mathbf{ c}} = E_\mathbf{d}(J(\mathbf{c}))$.
\end{proof}
\begin{cor}\label{cor} An orbit closure of type $A,D,E$ is smooth,
if and only if it is rationally smooth.
\end{cor}
\end{subsection}

\end{section}


\begin{section}{Geometric approach.}
We present in this section an independent proof for Theorem
\ref{smoothness}.
The interest of this proof is that we give an homological realization of the constants calculated in the previous section. To be more precise, the derivative of a coefficient $\Omega$ at $u=1$ can be seen as an
Euler-Poincar\'e cha\-rac\-te\-ris\-tic of a complex variety. The
principle on which the alternative proof works is that if a cone $X$ in $\mathbb{C}^n$
is rationally smooth, then, the rational cohomology of the projectivization
$\mathbb{P}(X)$ in
$\mathbb{P}^{n-1}$ is the same as the rational cohomology of
$\mathbb{P}^{d}$ where $d=\dim\mathbb{P}(X)$.

\begin{subsection}{Euler-Poincar\'e characteristic.}
In the following section, the varieties $\overline{\mathcal {O}}$ are considered on the field $\mathbb {C}$. We fix $\mathbf{ c}\in \mathbb{ N}^{\nu}$ of $\mathbf{ i}$-homogeneity
$\mathbf{ d}$. Set
$X_{\mathbf{ c}}:=\mathbb{P}
(\overline{\mathcal{O}_{\mathbf{c}}}\backslash\{0\})$,
as a complex variety and let
$X_{\mathbf{ c}}(\mathbb{F}_{q})$, resp. $X_{\mathbf{ c}}(\overline{\mathbb{F}_q})$, be the cor\-res\-ponding variety on the finite
field $\mathbb{F}_{q}$, resp $\overline{\mathbb{F}_q}$. We know that there exists a discrete valuation ring $R\subset\mathbb{C}$ with residue field of characteristic $p$ and a variety $\hat{X}_{\mathbf{ c}}$ defined on $R$ such that we get the variety $X_{\mathbf{ c}}(\overline{\mathbb{F}_q})$ over $\overline{\mathbb{F}_q}$ and the variety $X_{\mathbf{ c}}$ over $\mathbb{C}$ by the base change to Spec$\overline{\mathbb{F}_q}$ and Spec$\mathbb{C}$ respectively.\par
Let $\chi(X_{\mathbf{ c}})$ be the Euler-Poincar\'e characteristic of
the variety $X_{\mathbf{ c}}$. Then
\begin{theorem}\label{geometric} For every $\mathbf{ c}\in \mathbb{ N}^{\nu}$ of
$\mathbf{ i}$-homogeneity $\mathbf{ d}$ we have
$$\chi(X_{\mathbf{ c}})=\dim E_{\mathbf{d}}(J(\mathbf{c})).$$
\end{theorem}

\begin{proof}
We first transpose the calculation of $\chi(X_{\mathbf{ c}})$ in the
context of $l$-adic cohomo\-lo\-gy with compact support.
We have
$$\dim_{\mathbb{Q}}H^i(X_{\mathbf{ c}},\mathbb{Q})=
\dim_{\mathbb{C}}H^i(X_{\mathbf{ c}},\mathbb{C})=
\dim_{\mathbb{C}}H_c^i(X_{\mathbf{ c}},\mathbb{C})=
\dim_{\overline{\mathbb{Q}}_{\ell}}H_c^i
(X_{\mathbf{ c}}(\overline{\mathbb{F}_q}),\overline{\mathbb{Q}}_{\ell}),$$
where $l$ is prime to $p$.
Hence,
\begin{equation}\label{ladic}
\chi(X_{\mathbf{ c}})=\sum_i(-1)^i\dim H_c^i
(X_{\mathbf{ c}}(\overline{\mathbb{F}_q}),\overline{\mathbb{Q}}_{\ell}).\end{equation}
Fix a prime number $p$ such that $l$ is prime to $p$. We know that there is
an action of the Frobenius ${\rm Fr}$ on the $l$-adic cohomology. By the
Grothendieck trace formula,~\cite[7.10]{danilov1}, we know that
\begin{equation}\mid X_{\mathbf{ c}}({\mathbb{F}}_{p^e})\mid=\sum_{i=0}^{2
\dim X_{\mathbf{ c}}}
(-1)^i{\rm tr}({\rm Fr}^e,H_c^i
(X_{\mathbf{ c}}(\overline{\mathbb{F}_q}),\overline{\mathbb{Q}}_{\ell})).\end{equation}
Now, the theorem of Deligne,~\cite[8.21]{danilov1}, asserts that
$H_c^i(X_{\mathbf{ c}}(\overline{\mathbb{F}_q}),\overline{\mathbb{Q}}_{\ell})$ is filtered by
${\rm Fr}$-stable subspace $W_{i,j}$ such that the eigenvalues of the
Frobenius ${\rm Fr}$ on the successive quotients are powers of $p$ up to a
root of one. So, there exists an $N$ such that if $e$ is a multiple of $N$,
the eigenvalues of ${\rm Fr}^e$ are powers of $p^{e}$. For these $e$, we
have
\begin{equation}\mid X_{\mathbf{ c}}(\mathbb{F}_{p^e})\mid=\sum_{i,j}
(-1)^ip^{ek_{i,j}}\dim(W_{i,j}/W_{i,j-1}).\end{equation}
Now, we know that $\mid X_{\mathbf{ c}}(\mathbb{F}_{q})\mid$ is a polynomial
$P_{\mathbf{c}}(q)$. So, the previous equation provides a polynomial
equality which is true for an infinite number of $p^e$. \par
Hence,
$\mid X_{\mathbf{ c}}(\mathbb{F}_{q})\mid=\sum_{i,j}
(-1)^iq^{k_{i,j}}\dim(W_{i,j}/W_{i,j-1})$, which implies by
(\ref{ladic}) that $\chi(X_{\mathbf{ c}})=P_{\mathbf{c}}(1)$.
The proof of the theorem relies now on the following lemma.
\end{proof}

\begin{lemma}
For every $\mathbf{ c}\in \mathbb{ N}^{\nu}$ of $\mathbf{ i}$-homogeneity $\mathbf{ d}$ we
have
$$P_{\mathbf{c}}(1)=\dim(E_{\mathbf{d}}(J(\mathbf{c}))).$$
\end{lemma}

\begin{proof}
First of all, let $\mathbf{S}_{\mathbf{c}}$ be the set of elements $\mathbf{c'}$ of
$\mathbf{i}$-homogeneity $\mathbf{d}$ such that
there exists $\zS\in\sop$ such that
$\mathbf{c}^{min}+\mathbf{op}^\zS =\mathbf{c'}\preceq\mathbf{c}$. We have :
$$\mathbf{S}_{\mathbf{c}}=\{\mathbf{c}^{ij},\,(i,j)\in J(\mathbf{c})\},$$
with
$$c_s^{ij}=\left\{\begin{array}{ll}
1 & \textup{if } \alpha^s=
\alpha_i+\alpha_j\\
d_i-1 & \textup{if } \alpha^s=\alpha_i \\
d_j-1 & \textup{if } \alpha^s =
\alpha_j \\
0 & \textup{otherwise.}
\end{array}\right.
$$

We want to calculate $P_{\mathbf{c}}(q)$. As
$\overline{\mathcal{O}_{\mathbf{c}}}$ is a disjoint union of orbits
${\mathcal{O}_{\mathbf{c'}}}$ with $\mathbf{c'}\preceq\mathbf{c}$, we have
to count the cardinality of the set of $\mathbb{F}_q$-rational points
$\mathcal{O}_{\mathbf{c'}}(\mathbb{F}_q)$. This will be denoted by
$Q_{\mathbf{c'}}(q)$,
thus $P_\mathbf{c}(q)
={1 \over q-1}
\sum_{\mathbf{c}^{\textup{min}}\prec \mathbf{c}'\preceq\mathbf{c}}
Q_{\mathbf{c'}}(q)$.
Let $x'$ be a point of
$\mathcal{O}_{\mathbf{c'}}(\mathbb{F}_q)$ and
let $G_{\mathbf{d},x'}$ be the isotropy group of
$x'$ in $G_{\mathbf{d}}$. Then, $Q_{\mathbf{c'}}(q)={\mid
G_{\mathbf{d}}\mid\over\mid G_{\mathbf{d},x'}\mid}$ and it is
known that $G_{\mathbf{d},x'}$ is the group of automorphism of
the module $\mathbf{e}(\mathbf{c'})$. Hence, $\mid
G_{\mathbf{d},x'}\mid=a_{\mathbf{c}'}(q)$ is given by
(\ref{automorphism}). In particular, the multiplicity of $(q-1)$ in
$Q_{\mathbf{c'}}(q)$ is $\sum_{i=1}^n d_i-\sum_{s=1}^\nu c'_s$.

We find that every ${Q_{\mathbf{c'}}\over q-1}|_{q=1}$ is zero unless
$\mathbf{c'}\in\mathbf{S}_{\mathbf{c}}$.
Using the formula ${\mid Gl_{d_i}(\mathbb{F}_q)\mid\over\mid
Gl_{d_i-1}(\mathbb{F}_q)\mid}=(q^{d_i}-1)q^{d_i-1}$ and L'Hospital's rule,
we find
${Q_{\mathbf{c}^{ij}}\over q-1}|_{q=1}=d_id_j$.
By the decomposition $P_{\mathbf{c}}(q)={1\over q-1} \sum_{\mathbf{c^{{\rm
min}}}\prec\mathbf{c'}\preceq\mathbf{c}}Q_{\mathbf{c'}}(q)$, we obtain the
lemma.
\end{proof}
\end{subsection}
\begin{subsection}{Alternative proof.}
We now explain how the previous theorem gives a geometric version of the
algebraic proof given in~\cite{bedschi1} and in the present paper.
First, as $\overline{\mathcal{O}_{\mathbf{c}}}$ is a cone, we claim that the rational
smoothness of $\overline{\mathcal{O}_{\mathbf{c}}}$ implies that the
cohomology of $X_{\mathbf{c}}$ is the same as the cohomology of the projective space,
i.e. $\dim H^i(X_{\mathbf{c}},\mathbb{Q})=\dim H^i(\mathbb{P}^{
d({\mathbf{c}})-1},\mathbb{Q})$. Let's sketch a proof for this claim, fol\-lo\-wing \cite{brion1}.  By Remark \ref{bormac}, the rational smoothness property implies that $H^i(\overline{\mathcal{O}_{\mathbf{c}}},\overline{\mathcal{O}_{\mathbf{c}}}\backslash\{0\})=0$, $i\not=2d(\mathbf{c})$,  $H^{2d(\mathbf{c})}(\overline{\mathcal{O}_{\mathbf{c}}}, \overline{\mathcal{O}_{\mathbf{c}}}\backslash\{0\})=\mathbb{Q}$. As $\overline{\mathcal{O}_{\mathbf{c}}}$ is contractible, we obtain by a long exact sequence in relative cohomology, that $\overline{\mathcal{O}_{\mathbf{c}}}\backslash\{0\}$ is a rational cohomology sphere of dimension $2d(\mathbf{c})-1$. Now, let $S^1\subset\mathbb{C}^*$ acting naturally on the cone $\overline{\mathcal{O}_{\mathbf{c}}}$. It is clear that $X_{\mathbf{c}}$ and $(\overline{\mathcal{O}_{\mathbf{c}}}\backslash\{0\})/S^1$ have the same rational cohomology. Hence, a Gysin exact sequence gives
$H^0(X_{\mathbf{c}})\simeq H^2(X_{\mathbf{c}})\simeq\ldots H^{2(d(\mathbf{c})-1)}(X_{\mathbf{c}})\simeq\mathbb{Q}$ and $H^{2k+1}(X_{\mathbf{c}})=0$. This proves the claim.
\par
Now, $\chi(X_{\mathbf{c}})=\chi(\mathbb{P}^{d(\mathbf{c})-1})=d(\mathbf{c})$. By Theorem~\ref{geometric}, this implies that
$\dim\overline{O_{\mathbf{c}}}=\dim E_{\mathbf{d}}(J(\mathbf{c}))$,
and then $\overline{O_{\mathbf{c}}}=E_{\mathbf{d}}(J(\mathbf{c}))$. Once
again, this implies that $\overline{O_{\mathbf{c}}}$ is smooth.

\end{subsection}
\end{section}
\section*{Acknowledgements}
The first author wants to thank Michel Brion for suggesting the geometric approach and Markus Reineke for interesting and helpful email conversations.

The second author thanks Robert B\'edard for several discussions on the
subject.

\providecommand{\bysame}{\leavevmode\hbox to3em{\hrulefill}\thinspace}

\end{document}